# RESEARCH ANNOUNCEMENT

# A SPLITTING PROPERTY
# FOR SUBALGEBRAS OF TENSOR PRODUCTS

LIMING GE

ABSTRACT. We prove a basic result about tensor products of a $\text{II}_1$ factor with a finite von Neumann algebra and use it to answer, affirmatively, a question asked by S. Popa about maximal injective factors.

## 1. INTRODUCTION

The splitting property referred to in the title states that, in certain situations, a subalgebra $\mathcal{A}$ of the tensor product $\mathcal{A}_1 \otimes \mathcal{A}_2$ of two algebras $\mathcal{A}_1$ and $\mathcal{A}_2$ (with units) that contains $\mathcal{A}_1$ "splits" as a tensor product $\mathcal{A}_1 \otimes \mathcal{A}_0$, where $\mathcal{A}_0$ is a subalgebra of $\mathcal{A}_2$. If $\mathcal{A}_1$ is $M_n(\mathbb{C})$, the algebra of all $n \times n$ matrices over the complex numbers $\mathbb{C}$, and $\mathcal{A}_2$ is another complex matrix algebra, with the tensor product taken over $\mathbb{C}$, the splitting results from a simple algebraic calculation. If $\mathcal{A}_1$ and $\mathcal{A}_2$ are arbitrary, infinite-dimensional algebras over fields of arbitrary characteristic, with varying assumptions on their structure, the situation is less clear.

When $\mathcal{A}_1$ and $\mathcal{A}_2$ have topological and analytic structure and the tensor product $\mathcal{A}_1 \otimes \mathcal{A}_2$ is formed to reflect that structure, the splitting question becomes a deeper and more intricate one. The principal result of this note involves von Neumann algebras $\mathcal{R}_1$ and $\mathcal{R}_2$ and their von Neumann-algebra tensor product $\mathcal{R}_1 \bar{\otimes} \mathcal{R}_2$. Specifically, we prove the following result.

**Theorem A.** *If $\mathcal{M}$ is a factor of finite type, $\mathcal{R}$ is a finite von Neumann algebra, and $\mathcal{S}$ is a von Neumann subalgebra of $\mathcal{M} \bar{\otimes} \mathcal{R}$ that contains $\mathcal{M}$ ($= \mathcal{M} \bar{\otimes} \mathbb{C}I$), there is a von Neumann subalgebra $\mathcal{R}_0$ of $\mathcal{R}$ such that $\mathcal{S} = \mathcal{M} \bar{\otimes} \mathcal{R}_0$.*

In the preceding statement, we use "finite" in the sense of Murray and von Neumann [M-vN]: $\mathcal{M}$ is either a factor of type $\text{II}_1$ or of type $\text{I}_n$ (isomorphic to $M_n(\mathbb{C})$).

Using Theorem A, we answer, affirmatively, a question raised by Sorin Popa [P1]. Recall that an *injective* von Neumann algebra acting on a Hilbert space $\mathcal{H}$ is one that is the image of an idempotent, norm 1 mapping of $\mathcal{B}(\mathcal{H})$, the algebra of all bounded operators on $\mathcal{H}$. Popa asks [P1, Problem 4.5(1)]:

*If $\mathcal{M}_1$, $\mathcal{M}_2$ are type $II_1$ factors and $\mathcal{B}_1 \subset \mathcal{M}_1$, $\mathcal{B}_2 \subset \mathcal{M}_2$ are maximal injective von Neumann subalgebras, is $\mathcal{B}_1 \otimes \mathcal{B}_2$ maximal injective in $\mathcal{M}_1 \otimes \mathcal{M}_2$? Is this true at least for $\mathcal{M}_2 = \mathcal{B}_2 = \mathcal{R}$, where $\mathcal{R}$ is the hyperfinite $II_1$ factor?*

We show that the answer is affirmative in the case where $\mathcal{R}$ is the hyperfinite $II_1$ factor.

Further work is in progress, by the author and R. Kadison, on the splitting property when the restriction on the types of the von Neumann algebras is removed.

## 2. Methods and preliminary results

We make important use of the slice-map technique of Tomiyama [To] (especially as formulated in [K-R IV, Exercise 12.4.36]). Conditional expectation techniques and Schwartz projections [Sc] play a key role in the arguments. A description of these techniques follows.

From [K-R IV], let $\mathcal{R}$ and $\mathcal{S}$ be von Neumann algebras and $\rho$ and $\sigma$ be non-zero elements of $\mathcal{R}_\sharp$ and $\mathcal{S}_\sharp$, respectively. Then

(i) there is a unique element $\rho \otimes \sigma$ of $(\mathcal{R} \bar{\otimes} \mathcal{S})_\sharp$ such that $(\rho \otimes \sigma)(R \otimes S) = \rho(R)\sigma(S)$ for each $R$ in $\mathcal{R}$ and $S$ in $\mathcal{S}$ and that $\|\rho \otimes \sigma\| = \|\rho\|\|\sigma\|$;

(ii) there are unique operators $\Phi_\sigma(\tilde{T})$ and $\Psi_\rho(\tilde{T})$ in $\mathcal{R}$ and $\mathcal{S}$, respectively, corresponding to each $\tilde{T}$ in $\mathcal{R} \bar{\otimes} \mathcal{S}$, satisfying

$$\rho'(\Phi_\sigma(\tilde{T})) = (\rho' \otimes \sigma)(\tilde{T}), \qquad \sigma'(\Psi_\rho(\tilde{T})) = (\rho \otimes \sigma')(\tilde{T})$$

for each $\rho'$ in $\mathcal{R}_\sharp$ and each $\sigma'$ in $\mathcal{S}_\sharp$;

(iii) $\Phi_\sigma$ and $\Psi_\rho$ (as defined by (ii)) are ultraweakly continuous linear mappings of $\mathcal{R} \bar{\otimes} \mathcal{S}$ onto $\mathcal{R}$ and $\mathcal{S}$, respectively, satisfying

$$\Phi_\sigma((A \otimes I)\tilde{T}(B \otimes I)) = A\Phi_\sigma(\tilde{T})B,$$
$$\Psi_\rho((I \otimes C)\tilde{T}(I \otimes D)) = C\Psi_\rho(\tilde{T})D$$

for each $\tilde{T}$ in $\mathcal{R} \bar{\otimes} \mathcal{S}$; $A$, $B$ in $\mathcal{R}$; and $C$, $D$ in $\mathcal{S}$; and

$$\Phi_\sigma(R \otimes S) = \sigma(S)R, \qquad \Psi_\rho(R \otimes S) = \rho(R)S$$

when $R \in \mathcal{R}$ and $S \in \mathcal{S}$;

(iv) $\Phi_\sigma(\tilde{T}) \in \mathcal{R}_0$ and $\Psi_\rho(\tilde{T}) \in \mathcal{S}_0$ if $\tilde{T} \in \mathcal{R}_0 \bar{\otimes} \mathcal{S}_0$, where $\mathcal{R}_0$ and $\mathcal{S}_0$ are von Neumann subalgebras of $\mathcal{R}$ and $\mathcal{S}$, respectively;

(v) $\tilde{T} \in \mathcal{R}_0 \bar{\otimes} \mathcal{S}_0$ if $\Phi_{\sigma'}(\tilde{T}) \in \mathcal{R}_0$ and $\Psi_{\rho'}(\tilde{T}) \in \mathcal{S}_0$ for each $\sigma'$ in $\mathcal{S}_\sharp$ and each $\rho'$ in $\mathcal{R}_\sharp$;

(vi) $\Phi_\sigma(\tilde{T})$ appears, naturally, as a bounded linear functional on $\mathcal{R}_\sharp$ (thus, it is an element of $\mathcal{R}$). It is defined by

$$\Phi_\sigma(\tilde{T}): \rho' \mapsto (\rho' \otimes \sigma)(\tilde{T}) \qquad (\rho' \in \mathcal{R}_\sharp)$$

and

$$|(\Phi_\sigma(\tilde{T}))(\rho')| = |(\rho' \otimes \sigma)(\tilde{T})|$$
$$\leq \|\rho'\|\|\sigma\|\|\tilde{T}\|.$$

From this, (v) is valid under the assumption that $\sigma'$ and $\rho'$ may be any elements of families of normal states that generate, linearly, norm-dense subspaces of $\mathcal{S}$ and $\mathcal{R}$, respectively.

Suppose $\mathcal{R}$ is the ultraweak closure of an ascending union of finite-dimensional C*-algebras $\mathfrak{A}_n$ acting on a Hilbert space $\mathcal{H}$ (we call this union a *tower* for $\mathcal{R}$). Each $\mathfrak{A}_n$ has a finite group $\mathcal{U}_n$ of unitary elements that generates it linearly. In [Sc], Schwartz constructs a linear mapping of $\mathcal{B}(\mathcal{H})$ onto the commutant $\mathcal{R}'$ of $\mathcal{R}$ (a *Schwartz projection*) by averaging over the groups $\mathcal{U}_n$, successively, and passing to a Banach limit.

Given a tower $t$ for $\mathcal{R}$ and a free ultrafilter $p$, there is a norm 1, linear, idempotent mapping (a *conditional expectation*) $\Phi_{t,p}$ of $\mathcal{B}(\mathcal{H})$ onto $\mathcal{R}'$ such that

**Proposition B.** *Each $\Phi_{t,p}$ is a conditional expectation of $\mathcal{B}(\mathcal{H})$ onto $\mathcal{R}'$ proper on each von Neumann algebra $\mathcal{S}$ containing $\mathcal{R}$ (onto $\mathcal{S} \cap \mathcal{R}'$).*

For this proposition, we need

**Definition C.** A conditional expectation $\Psi$ of a von Neumann algebra $\mathcal{S}$ onto a subalgebra is *proper* when $\Psi(T) \in \mathrm{co}_\mathcal{S}(T)^-$ for each $T$ in $\mathcal{S}$, where $\mathrm{co}_\mathcal{S}(T)^-$ is the weak-operator closure of $\mathrm{co}_\mathcal{S}(T)$ and

$$\mathrm{co}_\mathcal{S}(T) = \{UTU^* : U \text{ a unitary in } \mathcal{S}\}.$$

In the next result, we describe a conditional expectation $\Phi$ of a finite von Neumann algebra $\mathcal{M}$ onto a von Neumann subalgebra $\mathcal{N}$ and describe some of its properties that will be needed.

**Theorem D.** *If $\mathcal{N}$ is a von Neumann subalgebra of a finite von Neumann algebra $\mathcal{M}$ with (normalized) trace $\tau$, there is a conditional expectation $\Phi$ of $\mathcal{M}$ onto $\mathcal{N}$, whose trace $\tau|\mathcal{N}$ we denote by "$\tau_0$", such that*

$$\tau(TA) = \tau_0(\Phi(T)A) \qquad (T \in \mathcal{M}, \ A \in \mathcal{N}).$$

*Moreover, $\Phi$ is the unique conditional expectation of $\mathcal{M}$ onto $\mathcal{N}$ such that*

$$(*) \qquad \tau(T) = \tau_0(\Phi(T)) \qquad (T \in \mathcal{M}).$$

*If $\mathcal{M}$ admits a proper conditional expectation $\Psi$ of $\mathcal{M}$ onto $\mathcal{N}$, then $\Psi = \Phi$.*

## 3. An outline of the proof

We use the notation of Theorem A and Section 2. We assume that $\mathcal{M}$ is a factor of type II$_1$. If $\mathcal{S}$ is a von Neumann subalgebra of $\mathcal{M} \bar{\otimes} \mathcal{R}$ containing $\mathcal{M} \bar{\otimes} \mathbb{C}I$ and $\mathcal{T} = \mathcal{M}' \cap \mathcal{S}$, what we want to show is that $\mathcal{M} \bar{\otimes} \mathcal{T} = \mathcal{S}$. Since $\mathcal{M} \bar{\otimes} \mathbb{C}I \subseteq \mathcal{S}$ and $\mathcal{T} \subseteq \mathcal{S}$, of course $\mathcal{M} \bar{\otimes} \mathcal{T} \subseteq \mathcal{S}$. Now we use the slice-map technique. If $T \in \mathcal{S}$, since $\Phi_\sigma$ maps onto $\mathcal{M}$, of course $\Phi_\sigma(T) \in \mathcal{M}$ for each $\sigma$. Define $\tau_H$ on $\mathcal{M}$ by $\tau_H(M) = \tau(HM)$ when $H, M \in \mathcal{M}$. We know that $\{\tau_H : H \in \mathcal{M}\}$ generate a norm-dense subspace of $\mathcal{M}_\sharp$ [K-R II, Theorem 7.2.3]. It remains to show that

$$\Psi_{\tau_H}(T) = \Psi_\tau((H \otimes I)T) \in \mathcal{T} \qquad (H \in \mathcal{M}).$$

But $\Psi_\tau$ lifts the trace from $\mathcal{R}$ to $\mathcal{M} \bar{\otimes} \mathcal{R}$, so $\Psi_\tau$ is the unique conditional expectation of $\mathcal{M} \bar{\otimes} \mathcal{R}$ onto $\mathcal{R}$ satisfying $(*)$ of Theorem D.

By a deep result of S. Popa (see [P2]), there is a hyperfinite subfactor $\mathcal{R}_0$ of $\mathcal{M}$ with trivial relative commutant. Since $\mathcal{R}_0$ is hyperfinite, each $\Phi_{t,p}$ (for $\mathcal{R}_0$) is proper, and the restriction of $\Phi_{t,p}$ to $\mathcal{M} \bar{\otimes} \mathcal{R}$ coincides with $\Psi_\tau$ (from Theorem D). Now $\mathrm{co}_{\mathcal{R}_0}((H \otimes I)T) \subseteq \mathcal{S}$—recall that $T \in \mathcal{S}$ and $H \otimes I \in \mathcal{M} \subseteq \mathcal{S}$. So

$$\Psi_{\tau_H}(T) = \Psi_\tau((H \otimes I)T) = \Phi_{t,p}((H \otimes I)T) \in \mathrm{co}_{\mathcal{R}_0}((H \otimes I)T)^- \subseteq \mathcal{S}.$$

But $\Phi_{t,p}$ maps onto $\mathcal{R}'_0$, so $\Phi_{t,p}((H \otimes I)T) \in \mathcal{R}'_0 \cap \mathcal{S} = \mathcal{M}' \cap \mathcal{S} = \mathcal{T}$ [K-R II, IV, 12.4.37]. This completes the argument.

The corollary that follows answers the Popa question noted in the introduction, even in a slightly generalized form. The proof of it follows directly from our Theorem A and properties of injective von Neumann algebras.

**Corollary E.** *Let $\mathcal{R}$ be the hyperfinite $\mathrm{II}_1$ factor and $\mathcal{S}$ a type $\mathrm{II}_1$ von Neumann algebra with a maximal injective von Neumann subalgebra $\mathcal{B}$. Then $\mathcal{R} \bar{\otimes} \mathcal{B}$ is maximal injective in $\mathcal{R} \bar{\otimes} \mathcal{S}$.*

## References

[K-R II, IV]  R. Kadison and J. Ringrose, *Fundamentals of the theory of operator algebras* II, Academic Press, Orlando, 1986; IV, Birkhäuser, Boston, 1992.

[M-vN]  F. Murray and J. von Neumann, *On rings of operators*, Ann. of Math. (2) **37** (1936), 116–229.

[P1]  S. Popa, *Maximal injective subalgebras in factors associated with free groups*, Adv. Math. **50** (1983), 27–48.

[P2]  ———, *On a problem of R. V. Kadison on maximal abelian *-subalgebras in factors*, Invent. Math. **65** (1981), 269–281.

[Sc]  J. Schwartz, *Two finite, non-hyperfinite, non-isomorphic factors*, Comm. Pure Appl. Math. **16** (1963), 19–26.

[To]  J. Tomiyama, *Tensor products and projections of norm one in von Neumann algebras*, Seminar notes, Math. Institute, University of Copenhagen, 1970.

Department of Mathematics, University of Pennsylvania, Philadelphia, Pennsylvania 19104-6395

*E-mail address*: liming@math.upenn.edu